\documentclass[11pt]{article}
\usepackage{amsmath}
\usepackage{amssymb}
\usepackage{amsfonts}
\usepackage{color}
\allowdisplaybreaks

\newtheorem{Theorem}{\textbf{Theorem}}[section]
\newtheorem{Lemma}{\textbf{Lemma}}[section]
\newtheorem{Proposition}{\textbf{Proposition}}[section]
\newtheorem{Corollary}{\textbf{Corollary}}[section]
\newtheorem{Remark}{\textbf{Remark}}[section]
\newtheorem{Example}{\textbf{Example}}[section]
\newtheorem{Definition}{\textbf{Definition}}[section]

\newenvironment{theorem}{\begin{Theorem}$\!\!\!$}{\end{Theorem}}
\newenvironment{lemma}{\begin{Lemma}$\!\!\!$}{\end{Lemma}}

\numberwithin{equation}{section}


\begin{document}
\title{Remarks on decay effects of regularity loss type wave equations with structural damping terms}
\author{Hironori Michihisa${}^\ast$\\ 
{\small Department of Mathematics, Graduate School of Science, Hiroshima University} \\
{\small Higashi-Hiroshima 739-8526, Japan}
}

\date{
}

\maketitle

\begin{abstract}
After GGH model was proposed in \cite{GGH}, 
Ikehata-Iyota \cite{IY} showed decay estimates for the total energy of solutions to GGH equations uniformly in the initial data. 
However, their results imply that the total energy is bounded when the initial data belong to the energy space. 
That is, whether it actually decays has not been known so far. 
In this paper we report a positive answer to that question.  
\end{abstract}

\footnote[0]{\hspace{-2em} ${}^\ast$Corresponding author.}
\footnote[0]{\hspace{-2em} \textit{Email:} hi.michihisa@gmail.com}
\footnote[0]{\hspace{-2em} 2010 \textit{Mathematics Subject Classification}. 35L05, 35L30, 35B65}
\footnote[0]{\hspace{-2em} \textit{Keywords and Phrases}: Damping, Regularity loss, Total energy, Diffusion structure}

\section{Introduction}
In this paper we consider the Cauchy problem of the solution to the wave equation with the structural damping term
\begin{equation}
\label{1.1}
\begin{cases}
u_{tt}-\Delta u+\Delta^2 u_t=0, & t>0,\quad x\in\textbf{R}^n, \\
u(0,x)=u_0(x),
\quad
u_t(0,x)=u_1(x),& x\in\textbf{R}^n, \\
\end{cases}
\end{equation}
where $n\ge1$ and $[u_0,u_1]\in (H^1(\textbf{R}^n)\cap L^1(\textbf{R}^n))\times(L^2(\textbf{R}^n)\cap L^1(\textbf{R}^n))$. 

Wave equations with structural damping terms $(-\Delta)^\theta u_t$ ($\theta\ge0$) are already studied by many mathematicians. 
For example, in Han-Milani \cite{HM}, Hosono \cite{H}, Hosono-Ogawa \cite{HO}, Karch \cite{K}, Marcati-Nishihara \cite{MN}, Matsumura \cite{M}, Michihisa \cite{Mi}, Narazaki \cite{Na}, Nishihara \cite{Ni}, Sakata-Wakasugi \cite{SW} and Takeda \cite{T}, they have investigated the diffusion structure for the case $\theta=0$, i.e., damped wave equations. 
Nishihara \cite{Ni} first studied the wave effect of the damped wave equation in the three dimensional case and higher order asymptotic expansions of the solution are found by Takeda \cite{T}, which inspired Michihisa \cite{Mi}. 
Karch \cite{K} obtained the asymptotic profile of the solution in the case $\theta\in[0,1/2)$ and his results are also applied to nonlinear problems in the same paper. 
See also D'Abbicco-Reissig \cite{DR} and references therein for recalling some decay estimates of solutions in the case $\theta\in(0,1]$. 

By the author's best knowledge, 
mathematical studies on the strongly damped wave equation ($\theta=1$) are started by Ponce \cite{P} and Shibata \cite{S}. 
After they obtained $L^p$-$L^q$ estimates of the solution, 
e.g., Ikehata-Natsume \cite{IN} showed the total energy decay estimates of the solution in the case $\theta\in[1/2,1]$. 
They also derived the $L^2$ decay estimates of the solution but the optimal estimates for the $L^2$ norm of the solution with $\theta=1$ were not obtained until Ikehata \cite{I} and Ikehata-Onodera \cite{IO}. 
Recently, the author \cite{Mi3} found the optimal leading term of the solution to the strongly damped wave equation in the Fourier space based on \cite{Mi2} where he  had established a method to appropriately expand the solution operators. 
Ikehata-Takeda \cite{IT} have classified asymptotic profiles of solutions to wave equations with damping terms $\mu (-\Delta)^\theta$  more precisely according to all the powers $\theta\in(0,1]$ and the constants $\mu>0$. 
See also Ikehata-Todorova-Yordanov \cite{ITY} where they stated existence results and asymptotic profiles of solutions to abstract equations including the strongly damped wave equation. 
These results above indicate that the solution of the strongly damped wave equation has a stronger wave-like property than that of the damped wave equation. 
The diffusion structure is, however, dominant in both cases. 

On the contrary, higher order fractional damping terms give a completely different phenomenon due to the regularity loss effect. 
After Ghisi-Gobbino-Haraux \cite{GGH} introduced the model equation \eqref{1.1}, 
Ikehata-Iyota \cite{IY} showed the decay estimates of the total energy and the optimal $L^2$ decay estimate of the solution. 
When they discussed the optimal decay estimate of the solution, they imposed some regularity conditions on the initial data to assure that the diffusive effect is still strong. 
They also considered the low regularity cases and found some decay estimates except for the case that the initial data belong to the energy space. 
In such a case, estimates derived by them tell us only that the $L^2$ norm of the solution is bounded. 
So we cannot determine whether it actually decays or not. 
In this paper, with simple calculations, we would like to give a positive answer to this question. 
\\

Before we state main results, 
we introduce some notation. 
The Fourier transform $\hat{f}(\xi)$ of a function $f(x)$ is defined by 
\[
\hat{f}(\xi)
:=(2\pi)^{-\frac{n}{2}}
\int_{\textbf{R}^n} 
e^{-ix\cdot\xi}f(x)\,dx.
\]
We write the $L^p(\textbf{R}^n)$ norm as $\|\cdot\|_p$ throughout this paper. 
When $p=2$, subscript $2$ is omitted.  
The total energy of the solution $u=u(t,x)$ to \eqref{1.1} is defined by 
\[
E(t):=\|u_t(t,\cdot)\|^2+\|\nabla u(t,\cdot)\|^2
=\|\hat{u}_t(t,\cdot)\|^2+\||\xi| \hat{u}(t,\cdot)\|^2.
\]
We also use Sobolev norms 
\[
\|D^\ell f\|
:=\||\xi|^\ell \hat{f}\|
\]
for $f\in H^\ell(\textbf{R}^n)$ with $\ell\ge0$. 
Now we are ready to state the main results of this paper. 

\begin{theorem}
\label{thm:1.1}
Let $n\ge1$ and $u$ be the solution to \eqref{1.1} with the initial data $[u_0,u_1]\in (H^1(\textbf{R}^n)\cap L^1(\textbf{R}^n))\times(L^2(\textbf{R}^n)\cap L^1(\textbf{R}^n))$. 
Then it holds that 
\begin{align*}
\lim_{t\to\infty}
E(t)
=0. 
\end{align*}
\end{theorem}

Theorem~\ref{thm:1.1} shows that the total energy does decay as $t\to\infty$. 
In this sense the term $(-\Delta)^2 u_t$ actually works as friction even for the low regularity case. 
However, it seems to be difficult to discuss the optimal decay rate. 
One may need another approach to obtain sharp decay estimates uniformly in the initial data. 

We also gain the following theorem first proved in \cite{IY} (see Theorem 1.1 in \cite{IY} with $\theta=2$) based on the proof of Theorem~\ref{thm:1.1} given in Section~\ref{sec:3}. 

\begin{theorem}
\label{thm:1.2}
Let $n\ge1$ and $\ell\ge0$. 
If $[u_0,u_1]\in (H^{\ell+1}(\textbf{R}^n)\cap L^1(\textbf{R}^n))\times(H^\ell(\textbf{R}^n)\cap L^1(\textbf{R}^n))$, 
then the solution $u$ to \eqref{1.1} satisfies 
\begin{align*}
E(t)
\le 
C(1+t)^{-\frac{n+2}{4}}\|u_0\|_1^2
+C(1+t)^{-\frac{n}{4}}\|u_1\|_1^2
+C(1+t)^{-\ell}(\|D^{\ell+1} u_0\|^2+\|D^\ell u_1\|^2), 
\quad t>0. 
\end{align*}
Here, $C>0$ is a constant independent of $t$ and the initial data. 
\end{theorem}
\section{Solution formula}
Ghisi-Gobbino-Haraux \cite{GGH} proved that problem \eqref{1.1} admits a unique mild solution in the class 
\[
u\in C([0,\infty); H^1(\textbf{R}^n))
\cap C^1([0,\infty); L^2(\textbf{R}^n))
\]
when the initial data belong to the energy space 
\[
[u_0,u_1]\in H^1(\textbf{R}^n)\times L^2(\textbf{R}^n).
\]
See also their results on smoothing properties. 
In this paper we will consider and fully use the formal solution formula in the Fourier space. 
The characteristic equation corresponding to problem \eqref{1.1} is 
\begin{align*}
\lambda^2+|\xi|^4 \lambda+|\xi|^2=0
\end{align*}
and we put its solutions as 
\begin{align*}
\lambda_\pm
:=\frac{-|\xi|^4\pm\sqrt{|\xi|^8-4|\xi|^2}}{2}.
\end{align*}
Hence the solution formula in the Fourier space is given by 
\begin{align*}
\hat{u}(t,\xi)
=E_0(t,\xi)\widehat{u_0}
+E_1(t,\xi)
\left(
\frac{|\xi|^4}{2}\widehat{u_0}
+\widehat{u_1}
\right),
\end{align*}
where 
\begin{align*}
E_0(t,\xi):=
\begin{cases}
\displaystyle{
e^{-\frac{1}{2}t|\xi|^4}
\cos
\left(
\frac{t\sqrt{4|\xi|^2-|\xi|^8}}{2}
\right)
},
& |\xi|\le\sqrt[3]{2}, \\[22pt]
\displaystyle{
e^{-\frac{1}{2}t|\xi|^4}
\cosh
\left(
\frac{t\sqrt{|\xi|^8-4|\xi|^2}}{2}
\right)
},
& |\xi|\ge\sqrt[3]{2}, 
\end{cases}
\end{align*}
\begin{align*}
E_1(t,\xi):=
\begin{cases}
\displaystyle{
e^{-\frac{1}{2}t|\xi|^4}
\sin
\left(
\frac{t\sqrt{4|\xi|^2-|\xi|^8}}{2}
\right)
\biggr/
\frac{\sqrt{4|\xi|^2-|\xi|^8}}{2}
},
& |\xi|\le\sqrt[3]{2}, \\[22pt]
\displaystyle{
e^{-\frac{1}{2}t|\xi|^4}
\sinh
\left(
\frac{t\sqrt{|\xi|^8-4|\xi|^2}}{2}
\right)
\biggr/
\frac{\sqrt{|\xi|^8-4|\xi|^2}}{2}
},
& |\xi|\ge\sqrt[3]{2}. 
\end{cases}
\end{align*}
Note that $\{\xi\in\textbf{R}^n: |\xi|=\sqrt[3]{2}\}$ is a removable singular set of $\hat{u}$ and so we can write $\hat{u}_t$ as follows:
\begin{equation}
\label{2.1}
\begin{split}
& \hat{u}_t(t,\xi) \\
& =-2
\frac{|\xi|}{\sqrt{4-|\xi|^6}}
\sin
\left(
\frac{t\sqrt{4|\xi|^2-|\xi|^8}}{2}
\right)
e^{-\frac{1}{2}t|\xi|^4}
\widehat{u_0} \\
& \quad 
+\left[
-\frac{|\xi|^3}{\sqrt{4-|\xi|^6}} 
\sin
\left(
\frac{t\sqrt{4|\xi|^2-|\xi|^8}}{2}
\right)
+\cos
\left(
\frac{t\sqrt{4|\xi|^2-|\xi|^8}}{2}
\right)
\right]
e^{-\frac{1}{2}t|\xi|^4}
\widehat{u_1},
\qquad 
|\xi|\le\sqrt[3]{2},
\end{split}
\end{equation}
\begin{equation}
\label{2.2}
\begin{split}
& \hat{u}_t(t,\xi) \\
& =-\frac{2}{|\xi|^2\sqrt{1-4/|\xi|^6}}
\sinh
\left(
\frac{t\sqrt{|\xi|^8-4|\xi|^2}}{2}
\right)
e^{-\frac{1}{2}t|\xi|^4}
\widehat{u_0} \\
& \quad
+\left[
-\frac{1}{\sqrt{1-4/|\xi|^6}}
\sinh
\left(
\frac{t\sqrt{|\xi|^8-4|\xi|^2}}{2}
\right)
+\cosh
\left(
\frac{t\sqrt{|\xi|^8-4|\xi|^2}}{2}
\right)
\right]
e^{-\frac{1}{2}t|\xi|^4}
\widehat{u_1}, 
\qquad 
|\xi|\ge\sqrt[3]{2}.
\end{split}
\end{equation}

\section{Proofs}
\label{sec:3}
First, we prepare some basic estimates. 
Since they can be shown by simple calculations, proofs are omitted. 
\begin{lemma}
\label{lem:3.1}
Let $n\ge1$ and $m\ge0$. 
Then there exists a constant $C>0$ such that 
\begin{align}
\label{3.1}
\||\xi|^m e^{-\frac{1}{2}t|\xi|^4}\|_{L^2(|\xi|\le1)}
\le C(1+t)^{-\frac{n+2m}{8}},
\quad 
t\ge0,
\end{align}
\begin{align}
\label{3.2}
\sup_{|\xi|\ge1}
\frac{e^{-\frac{t}{|\xi|^2}}}{|\xi|^m}
\le C(1+t)^{-\frac{m}{2}},
\quad 
t\ge0.
\end{align}
\end{lemma}

In the low frequency region one can observe a diffusion wave property. 
Due to its diffusion structure the following decay estimate can be easily proved. 
\begin{lemma}
\label{lem:3.2}
Let $n\ge1$ and $u$ be the solution to \eqref{1.1} with the initial data $[u_0,u_1]\in (H^1(\textbf{R}^n)\cap L^1(\textbf{R}^n))\times(L^2(\textbf{R}^n)\cap L^1(\textbf{R}^n))$. 
Then it holds that 
\begin{align*}
\|\hat{u}_t(t,\cdot)\|_{L^2(|\xi|\le1)}^2
+\||\xi|\hat{u}(t,\cdot)\|_{L^2(|\xi|\le1)}^2
\le C(1+t)^{-\frac{n}{4}}\|u_0\|_1^2 
+C(1+t)^{-\frac{n+2}{4}}\|u_1\|_1^2, 
\quad t>0.
\end{align*}
Here, $C>0$ is a constant independent of $t$ and the initial data. 
\end{lemma}
\textbf{Proof.} 
From \eqref{2.1}, we see that 
\[
|\hat{u}_t(t,\xi)|
\le C|\xi| e^{-\frac{1}{2}t|\xi|^4}\|u_0\|_1
+Ce^{-\frac{1}{2}t|\xi|^4}\|u_1\|_1
\]
for $t>0$ and $\xi\in\textbf{R}^n$ with $|\xi|\le1$. 
It follows from \eqref{3.1} that 
\begin{align*}
\|\hat{u}_t(t,\cdot)\|_{L^2(|\xi|\le1)}^2
\le  C(1+t)^{-\frac{n+2}{4}}\|u_0\|_1^2
+C(1+t)^{-\frac{n}{4}}\|u_1\|_1^2, 
\quad t>0.
\end{align*}
Similarly, we have 
\[
|\xi||\hat{u}(t,\xi)|
\le C|\xi| e^{-\frac{1}{2}t|\xi|^4}\|u_0\|_1
+Ce^{-\frac{1}{2}t|\xi|^4}\|u_1\|_1
\]
for $t>0$ and $\xi\in\textbf{R}^n$ with $|\xi|\le1$. 
Thus, as for $\nabla u$, we obtain the same $L^2$ estimate as above. 
Therefore the proof is complete. 
$\Box$
\\

Next lemma says that the total energy in the middle frequency region decays exponentially in $t$. 

\begin{lemma}
\label{lem:3.3}
Let $n\ge1$ and $u$ be the solution to \eqref{1.1} with the initial data $[u_0,u_1]\in H^1(\textbf{R}^n)\times L^2(\textbf{R}^n)$. 
Then it holds that 
\begin{align*}
\|\hat{u}_t(t,\cdot)\|_{L^2(1\le|\xi|\le\sqrt{2})}^2
+\||\xi|\hat{u}(t,\cdot)\|_{L^2(1\le|\xi|\le\sqrt{2})}^2
\le Ce^{-ct}\|u_0\|_2^2 
+Ce^{-ct}\|u_1\|_2^2, 
\quad t>0.
\end{align*}
Here, $C>0$ and $c>0$ are constants independent of $t$ and the initial data. 
\end{lemma}
\textbf{Proof.} 
To derive the desired estimate we can use the estimates below: 
\[
\mbox{(i)}\,\,\,
\sup_{x>0}
\left|
\frac{\sin x}{x}
\right|
=1,
\qquad
\mbox{(ii)}\,\,\,
\frac{\sinh x}{x}\le Ce^x 
\quad (x>0),
\]
\[
\mbox{(iii)}\,\,\,
0<
e^{-\frac{1}{2}t|\xi|^4}
\sinh
\left(
\frac{t\sqrt{|\xi|^8-4|\xi|^2}}{2}
\right)
\le e^{\lambda_+ t},
\quad
0<
e^{-\frac{1}{2}t|\xi|^4}
\cosh
\left(
\frac{t\sqrt{|\xi|^8-4|\xi|^2}}{2}
\right)
\le e^{\lambda_+ t}.
\]
From \eqref{2.1} with (i), one has 
\[
|\hat{u}_t(t,\xi)|
\le C(1+t)e^{-\frac{t}{2}}
(|\widehat{u_0}(\xi)|+|\widehat{u_1}(\xi)|)
\]
for $t>0$ and $\xi\in\textbf{R}^n$ with $1\le|\xi|\le\sqrt[3]{2}$. 
Similarly, it follows from \eqref{2.2} with (ii) and (iii) that 
\[
|\hat{u}_t(t,\xi)|
\le C(1+t) e^{\lambda_+ t}
(|\widehat{u_0}(\xi)|+|\widehat{u_1}(\xi)|) 
\]
for $t>0$ and $\xi\in\textbf{R}^n$ with $\sqrt[3]{2}\le|\xi|\le4$. 
Note that 
\[
\lambda_+
\le -\frac{2}{|\xi|^2}
\le -\frac{1}{8}
\]
if $\sqrt[3]{2}\le|\xi|\le4$. 
So, for arbitrary $c\in(0,1/8)$, there exists a constant $C>0$ such that 
\[
|\hat{u}_t(t,\xi)|
\le 
Ce^{-ct}
(|\widehat{u_0}(\xi)|+|\widehat{u_1}(\xi)|) 
\]
for $t>0$ and $\xi\in\textbf{R}^n$ with $1\le|\xi|\le4$. 
The same estimate also holds for $|\xi||\hat{u}(t,\xi)|$ and thus the lemma is obtained. 
$\Box$
\\

The following lemma is the essential part of this paper. 
Especially, the behavior of $\lambda_+$ in the high frequency region is a key to derive desired decay estimates for the total energy. 
See \eqref{3.5} below which comes from the negativity of $\lambda_+$. 
\begin{lemma}
\label{lem:3.4}
Let $n\ge1$ and $u$ be the solution to \eqref{1.1} with the initial data $[u_0,u_1]\in H^1(\textbf{R}^n)\times L^2(\textbf{R}^n)$. 
Then it holds that 
\begin{align*}
\lim_{t\to\infty}
\biggr(
\|\hat{u}_t(t,\cdot)\|_{L^2(|\xi|\ge1)}^2
+\||\xi|\hat{u}(t,\cdot)\|_{L^2(|\xi|\ge1)}^2
\biggr)
=0.
\end{align*}
\end{lemma}
\textbf{Proof.} 
First, we see that 
\begin{align*}
\lambda_+
=-\frac{1}{|\xi|^2}\frac{2}{1+\sqrt{1-4/|\xi|^6}}
\le -\frac{1}{|\xi|^2} 
\end{align*}
for $t>0$ and $\xi\in\textbf{R}^n$ with $|\xi|\ge\sqrt[3]{2}$. 
Thus it follows from \eqref{2.2} and (iii) in the previous proof that 
\begin{align}
\label{3.3}
|\hat{u}_t(t,\xi)|
\le C \frac{e^{-\frac{t}{|\xi|^2}}}{|\xi|^2} 
|\widehat{u_0}(\xi)|
+Ce^{-\frac{t}{|\xi|^2}}|\widehat{u_1}(\xi)|
\end{align}
for $t>0$ and $\xi\in\textbf{R}^n$ with $|\xi|\ge\sqrt{2}$. 
Similarly, one also has 
\begin{align}
\label{3.4}
|\xi||\hat{u}(t,\xi)|
\le C e^{-\frac{t}{|\xi|^2}} 
|\xi| |\widehat{u_0}(\xi)|
+C \frac{e^{-\frac{t}{|\xi|^2}}}{|\xi|^3}|\widehat{u_1}(\xi)|
\end{align}
for $t>0$ and $\xi\in\textbf{R}^n$ with $|\xi|\ge\sqrt{2}$. 
For fixed $\xi\in\textbf{R}^n$ with $|\xi|\ge\sqrt{2}$, we can easily see that 
\begin{align}
\label{3.5}
0<e^{-2\frac{t}{|\xi|^2}}\le1, 
\qquad
\lim_{t\to\infty}e^{-2\frac{t}{|\xi|^2}}=0. 
\end{align}
Therefore the Lebesgue dominated convergence theorem gives the desired result under the condition that the initial data belong to the energy space. 
$\Box$
\\

Theorem~\ref{thm:1.1} is a direct result of Lemmas~\ref{lem:3.2}-\ref{lem:3.4}. 
At the end of this paper, we prove Theorem~\ref{thm:1.2} based on the proof of Theorem~\ref{thm:1.1}. 
\\

\textbf{Proof of Theorem~\ref{thm:1.2}.} 
Recall inequalities \eqref{3.3} and \eqref{3.4} to derive 
\[
|\hat{u}_t(t,\xi)|
\le C \frac{e^{-\frac{t}{|\xi|^2}}}{|\xi|^{\ell+3}} 
|\xi|^{\ell+1} |\widehat{u_0}(\xi)|
+C\frac{e^{-\frac{t}{|\xi|^2}}}{|\xi|^\ell}
|\xi|^\ell |\widehat{u_1}(\xi)|,
\]
\[
|\xi||\hat{u}(t,\xi)|
\le C \frac{e^{-\frac{t}{|\xi|^2}}}{|\xi|^\ell} 
|\xi|^{\ell+1} |\widehat{u_0}(\xi)|
+C \frac{e^{-\frac{t}{|\xi|^2}}}{|\xi|^{\ell+3}}
|\xi|^\ell |\widehat{u_1}(\xi)|,
\]
for $t>0$ and $\xi\in\textbf{R}^n$ with $|\xi|\ge\sqrt{2}$. 
With the aid of \eqref{3.2}, we have 
\[
\|\hat{u}_t(t,\cdot)\|_{L^2(|\xi|\ge\sqrt{2})}^2 
\le C(1+t)^{-(\ell+3)}\|D^{\ell+1} u_0\|^2 
+C(1+t)^{-\ell}\|D^\ell u_1\|^2, 
\]
\[
\||\xi|\hat{u}(t,\cdot)\|_{L^2(|\xi|\ge\sqrt{2})}^2
\le C(1+t)^{-\ell}\|D^{\ell+1} u_0\|^2 
+C(1+t)^{-(\ell+3)}\|D^\ell u_1\|^2, 
\]
for $t>0$. 
Combining these estimates with Lemmas~\ref{lem:3.2} and \ref{lem:3.3}, 
we complete the proof of the theorem. 
$\Box$ 
\\

\noindent 
\textbf{Acknowledgement.} 
The author would like to thank Professor Ryo Ikehata for useful suggestions and continuing support. 

\bibliographystyle{amsplain}

\end{document}